\documentclass[12pt]{article}
\usepackage{amssymb}
\usepackage{amsthm}
\usepackage{amsmath}
\usepackage{amsfonts}

\begin{document}

\theoremstyle{definition}
\newtheorem{defn}{Definition}[section]
\newtheorem{conv}[defn]{Convention}
\newtheorem{rmk}[defn]{Remark}
\theoremstyle{plain}
\newtheorem{thm}[defn]{Theorem} 
\newtheorem{cor}[defn]{Corollary}
\newtheorem{lemma}[defn]{Lemma}
\newtheorem{propn}[defn]{Proposition}
\newtheorem{add}[defn]{Addendum}

\newcommand{\rxm}{\ensuremath{\mathbb{R} \times M}}
\newcommand{\reals}{\ensuremath{\mathbb{R}}}
\newcommand{\half}{\ensuremath{\frac{1}{2}}}
\newcommand{\into}{\ensuremath{\hookrightarrow}}
\newcommand{\tb}{\mathop{\rm tb}\nolimits}
\newcommand{\interior}{\mathop{\rm interior}\nolimits}
\newcommand{\closure}{\mathop{\rm closure}\nolimits}
\newcommand{\id}{\mathop{\rm id}\nolimits}
\newcommand{\supp}{\mathop{\rm supp}\nolimits}
\newcommand{\Span}{\mathop{\rm span}\nolimits}
\newcommand{\pf}{\mathop{\rm pf}\nolimits}

\title{Explicit concave fillings of contact three-manifolds}

\author{David T. Gay\\
Mathematics Department, University of Arizona\\
617 N. Santa Rita, P.O. Box 210089\\
Tucson AZ 85721, U.S.A.\\
email: dtgay@math.arizona.edu
}

%\begin{abstract}

%When $(M,\xi)$ is a contact $3$-manifold we say that a compact
%symplectic $4$-manifold $(X,\omega)$ is a concave filling of $(M,\xi)$
%if $M = -\partial X$ and if there exists a Liouville vector field $V$
%defined on a neighborhood of $M$, transverse to $M$ and pointing {\em
%in} to $X$, such that $\xi$ is the kernel of $\imath_V \omega$
%restricted to $M$. We give explicit, handle-by-handle constructions of
%concave fillings of all closed, oriented, contact $3$-manifolds.
%\end{abstract}

\maketitle
%\begin{spacing}{2}
\section{Introduction}

Suppose that $(X,\omega)$ is a symplectic manifold and that there
exists a Liouville vector field $V$ defined in a neighborhood of and
transverse to $M = \partial X$. Then $V$ induces a contact form
$\alpha = \imath_V \omega|_M$ on $M$ which determines the germ of
$\omega$ along $M$. One should think of the contact manifold
$(M,\xi=\ker \alpha)$ as controlling the behavior of $\omega$ ``at
infinity''. If $V$ points out of $X$ along $M$ then we call
$(X,\omega)$ a {\em convex filling} of $(M,\xi)$, and if $V$ points
into $X$ along $M$ then we call $(X,\omega)$ a {\em concave filling}
of $(M,\xi)$.

Much attention has been given in recent years to constructions of
convex fillings of contact $3$-manifolds; see~\cite{GompfStein} for
example.  Etnyre and Honda~\cite{EtnyreHonda} very recently began a
careful investigation of concave fillings and proved that every
compact, oriented, contact $3$-manifold has (infinitely many) concave
fillings, but their proof depends on a result of Lisca and
Matic~\cite{LiscaMatic} that every Stein surface embeds as a domain in
a closed K\"ahler manifold. This result is not explicitly constructive
and in particular does not give handlebody decompositions of the
concave fillings. In this paper we present a method to explicitly
construct, handle by handle, concave fillings of contact
$3$-manifolds, without reference to the result of Lisca and Matic. The
main theorem we prove is:
\begin{thm} \label{T:ConcaveFillings}
Every closed, oriented, contact $3$-manifold has a concave filling
and, furthermore, a handlebody decomposition of the filling can be
given explicitly in terms of the contact structure.
\end{thm}

While the author was writing this paper Akbulut and
Ozbagci~\cite{AkbOzb} presented, using different techniques, a
constructive proof of the fact that every Stein surface embeds in a
closed symplectic $4$-manifold. The reader will also find some common
themes in an earlier paper by Akbulut and Ozbagci~\cite{AkbOzb2}.

The author would like to thank Emmanuel Giroux for helpful
correspondence regarding his recent results, Ichiro Torisu for
pointing out Giroux's results to the author and the Nankai Institute
of Mathematics for support and hospitality.

\section{Tools and building blocks}

Henceforth in this paper, unless explicitly stated otherwise, we adopt
the following conventions: All manifolds are compact and oriented, and
all $3$-manifolds are closed. All contact structures are positive, and
all symplectic manifolds are oriented by their symplectic
structures. For the basic definitions in symplectic and contact
topology the reader is referred to~\cite{McDSal}
and~\cite{ThomasEliashbergGiroux}.

We will build our concave fillings by glueing together symplectic
cobordisms, which we now define. Let $(X,\omega)$ be a symplectic
$4$-manifold and let $M_1$ and $M_2$ be two $3$-manifolds with
$\partial X = (-M_1) \amalg M_2$. Suppose there exists a Liouville
vector field $V$ defined on a neighborhood of and transverse to
$\partial X$, pointing in along $M_1$ and out along $M_2$, and let
$\xi_i$ be the induced contact structure on $M_i$. In this situation,
following~\cite{EliashGivenHofer} and~\cite{EtnyreHonda}, we call
$(X,\omega)$ a symplectic cobordism from $(M_1,\xi_1)$ to
$(M_2,\xi_2)$ and we indicate the existence of such a cobordism with
the notation $(M_1,\xi_1) \prec (M_2,\xi_2)$. Note that this relation
is not reflexive, but it is transitive: If $(X_i,\omega_i)$ is a
symplectic cobordism from $(M_i,\xi_i)$ to $(M_{i+1},\xi_{i+1})$, for
$i\in\{1,2\}$, then it is possible to glue $(X_1,\omega_1)$ to
$(X_2,\omega_2)$ (after attaching a symplectic collar on
$(X_1,\omega_1)$ and possibly multiplying $\omega_2$ by a positive
constant) to form a symplectic cobordism from $(M_1,\xi_1)$ to
$(M_3,\xi_3)$.

A convex filling of $(M,\xi)$ is thus a symplectic cobordism from
$\emptyset$ to $(M,\xi)$ while a concave filling of $(M,\xi)$ is a
symplectic cobordism from $(M,\xi)$ to $\emptyset$.

Given any contact $3$-manifold $(M,\xi)$, the most basic building
block is a cobordism $(X,\omega)$ from $(M,\xi)$ to itself constructed
as follows: Let $\alpha$ be any contact form for $\xi$ and let $\omega
= d(e^s \alpha)$ on $\reals \times M$, where $s$ is the $\reals$
coordinate. Then $\omega$ is a symplectic form and $\partial_s$ is a
Liouville vector field inducing the contact form $e^f \alpha$ on the
graph of any function $f\colon\thinspace M \to \reals$, so for any two
functions $f < g$ on $M$ we can construct a cobordism $X = \{(s,p)
\mid f(p) \leq s \leq g(p)\} \subset \reals \times M$. Use of this
building block will often be assumed without explicit mention, but it
is essential for all of our constructions, and in particular for the
fact that cobordisms can be glued together.

Our main construction depends on the recent discovery of close
connections between contact structures and open book decompositions,
culminating in the work of Giroux~\cite{Giroux}. To discuss open book
decompositions, we begin with conventions regarding mapping class
groups. If $\Sigma$ is a compact surface with boundary, we take the
mapping class group of $\Sigma$ to be the group of
orientation-preserving self-homeomorphisms of $\Sigma$ fixing
$\partial \Sigma$ pointwise, modulo isotopies fixing $\partial \Sigma$
pointwise. We denote this mapping class group by
$\mathcal{M}(\Sigma)$. We multiply elements of $\mathcal{M}(\Sigma)$
in the same order in which we compose functions, and in fact we will
generally blur the distinction between a homeomorphism
$h\colon\thinspace \Sigma \to \Sigma$ and its equivalence class $h \in
\mathcal{M}(\Sigma)$. If $\Sigma^\prime \supset \Sigma$ is another
compact surface, there is a natural inclusion $\mathcal{M}(\Sigma)
\subset \mathcal{M}(\Sigma^\prime)$ given by extending $h \in
\mathcal{M}(\Sigma)$ to be the identity on $\Sigma^\prime \setminus
\Sigma$.

By ``polar coordinates'' on a solid torus $D^2 \times S^1$ we will
mean coordinates $(r,\mu,\lambda)$, where $(r,\mu)$ are polar
coordinates on $D^2$ and $\lambda \in S^1$. By polar coordinates near
a link $L$ we will mean polar coordinates on a neighborhood of each
component of $L$ such that $L = \{r=0\}$.

An open book decomposition of a $3$-manifold $M$ is a non-empty link
$L \subset M$ and a fibration $p\colon\thinspace M \setminus L \to
S^1$ such that near $L$ we can find polar coordinates
$(r,\mu,\lambda)$ with respect to which $p = \mu$. The link $L$ is
called the ``binding'' and the union of a fiber with $L$ is called a
``page''.  From such a structure we get the following data:
\begin{enumerate}
\item The topological type of a page, a compact, oriented surface
$\Sigma$ with $\partial \Sigma = L$.
\item The monodromy $h$, the element of $\mathcal{M}(\Sigma)$
represented by the return map for flow along any vector field $V$ on
$M \setminus L$ which is transverse to the pages and is meridinal near
$L$ (i.e. $V = \partial_\mu$ in polar coordinates near $L$).
\end{enumerate}
We will often refer to the pair $(\Sigma,h)$ as the open book
decomposition, and the manifold with this open book decomposition as
$B{(\Sigma,h)}$. We identify $B{(\Sigma,h)}$ with $(\Sigma \times
[0,2\pi])/\!\sim$ where $(p,2\pi) \sim (h(p),0)$ for all $p \in \Sigma$
and $(p,s) \sim (p,t)$ for all $p \in \partial \Sigma$ and all $s,t
\in [0,2\pi]$. The pages are $\Sigma \times \{t\} / \!\sim$ and the
binding is $L = (\partial \Sigma \times \{t\}) / \!\sim$. Notice that
$B{(\Sigma,h)} = B{(-\Sigma,h^{-1})} = - B{(\Sigma,h^{-1})}$.

A contact vector field for a contact structure $\xi$ is a vector field
$V$ such that flow along $V$ preserves $\xi$.  A Reeb vector field for
$\xi$ is a contact vector field which is transverse to $\xi$. Every
Reeb vector field $V$ arises from a contact form $\alpha$ for $\xi$ as
the unique vector field satisfying $\imath_V d\alpha =0$ and
$\alpha(V) = 1$.
\begin{defn} \label{D:Supports}
Let us say that a vector field $V$ on $M$ respects a given open book
decomposition of $M$ if $V$ is transverse to the pages, tangent to the
binding, and satisfies the following orientation condition: $V$
co-orients and hence orients each page, and also orients the binding,
and we require that this orientation of the binding agree with its
orientation as the boundary of a page.  A contact structure $\xi$ on
$M$ is {\em supported} by a given open book decomposition if there
exists a Reeb vector field for $\xi$ which respects the open book
decomposition
\end{defn} 
The first result below is due, in various versions, to
Giroux~\cite{Giroux}, Torisu~\cite{Torisu}, and Thurston and
Winkelnkemper~\cite{ThurstonWinkel}
\begin{thm}[]
\label{T:OBDGivesCS} 
Every open book decomposition of a $3$-manifold $M$ supports a unique
(up to isotopy) contact structure on $M$.
\end{thm}
Thurston and Winkelnkemper proved existence without mentioning
uniqueness. Torisu proved existence and uniqueness, but using a
different relationship between the contact structure and the open book
decomposition defined in terms of convex Heegard splittings. We will
discuss Torisu's version in section~\ref{S:ConstructBlocks}. Finally
Giroux pointed out a simple argument for the uniqueness. This result
becomes most powerful when coupled with the following theorem.
\begin{thm}[Giroux~\cite{Giroux}] \label{T:CSgivesOBD}
Every contact structure on a $3$-manifold $M$ is supported by an open
book decomposition of $M$.
\end{thm}

Given an open book decomposition $(\Sigma,h)$, we will refer to the
unique supported contact structure on $B(\Sigma,h)$ as
$\xi(\Sigma,h)$, and we will abbreviate with the notation
$\mathfrak{B}(\Sigma,h) = (B(\Sigma,h),\xi(\Sigma,h))$. 

\begin{rmk} \label{R:SThree}
If $\Sigma$ is a disk and $h=1$ then $B{(\Sigma,h)} = S^3$ and
$\xi{(\Sigma,h)}$ is the standard contact structure $\xi_0$ on
$S^3$. The $4$-ball with its standard symplectic form is a symplectic
cobordism from $\emptyset$ to $(S^3,\xi_0)$.
\end{rmk}

\begin{rmk} \label{R:SurfaceCrossDisk}
As a generalization, if $\Sigma$ is any compact surface with $\partial
\Sigma \neq \emptyset$ and $h=1$, and if we let $X$ be $\Sigma \times
D^2$ with ``rounded corners'' and $\omega$ be the sum of a volume form
on $\Sigma$ and a volume form on $D^2$, then $(X,\omega)$ is a
symplectic cobordism form $\emptyset$ to $\mathfrak{B}{(\Sigma,h)}$. We
will call this the standard convex filling of
$\mathfrak{B}{(\Sigma,1)}$.
\end{rmk}

The strategy of this paper is to translate cobordism questions for
contact $3$-mani\-folds into questions about the relationship between
mapping class group elements of various surfaces. We will rely on the
following fundamental fact:
\begin{thm}[Dehn~\cite{Dehn} and Lickorish~\cite{Lickorish}]
\label{T:Dehn} 
The mapping class group of a surface is generated by Dehn twists.
\end{thm}
It is important to distinguish between right-handed and left-handed
Dehn twists. A right-handed Dehn twist about a simple closed curve $C$
is the twist $\tau_C$ such that, if $\gamma$ is an arc transverse to
$C$ and we travel along $\gamma$ towards $C$, $\tau_C(\gamma)$
diverges from $\gamma$ by forking off to the right, going around $C$
and then rejoining $\gamma$ on the other side of $C$. This distinction
depends on the orientation of $\Sigma$ but not of $C$. We will shorten
the phrase ``right-handed Dehn twist'' to ``right twist''. If we are
given a single curve $C$, we will use the notation $\tau_C$ to refer
to the right twist about $C$. If we are given a sequence of curves
$\{C_i\}$, we will use the notation $\tau_i$ to refer to the right
twist about $C_i$. It will also be important to distinguish between
twists about homologically trivial and nontrivial curves. We will call
a twist about a homologically nontrivial curve a homologically
nontrivial twist.

We now present the fundamental building blocks for our
constructions. Throughout, suppose we are given a pair $(\Sigma,h)$,
where $\Sigma$ is a compact surface with $\partial \Sigma \neq
\emptyset$ and $h \in \mathcal{M}(\Sigma)$.

The first two results below are more or less immediate consequences of
results of Eliashberg~\cite{EliashStein} and Weinstein~\cite{Weinstein}
on symplectic handlebodies, reinterpreted in the context of supporting
open book decompositions. As stated here the results are probably well
known to experts, but for completeness we will provide proofs in
section~\ref{S:ConstructBlocks}.
\begin{propn} \label{P:OneHandle}
Let $\Sigma^\prime$ be $\Sigma$ with a $2$-dimensional $1$-handle
attached (at points $p,q \in \partial \Sigma$). Then
$\mathfrak{B}{(\Sigma,h)} \prec \mathfrak{B}{(\Sigma^\prime,h)}$, and
the cobordism is diffeomorphic to $[-1,0] \times B{(\Sigma,h)}$ with a
$4$-dimensional $1$-handle attached to $\{0\} \times B{(\Sigma,h)}$
(at the corresponding points $p$ and $q$ in the binding).
\end{propn}
Note that this result is straightforward if we ignore the contact and
symplectic structures; simply observe that the pages intersect the
boundary $S^2$ of each foot of the $1$-handle in longitudes running
from the north to the south pole, and that these extend across the
surgery as $2$-dimensional $1$-handles attached to each page

To discuss $2$-handles, notice that a knot $K$ lying in a page of an open
book decomposition is given a framing by a vector tangent to the page
and transverse to the knot; call this the page-framing of the knot and
abbreviate it $\pf(K)$.
\begin{propn} \label{P:NegTwoHandle}
Let $h^\prime = h \circ \tau_C$, where $\tau_C$ is a right twist along
a homologically nontrivial curve $C \subset \Sigma \setminus \partial
\Sigma$. Then $\mathfrak{B}{(\Sigma,h)} \prec
\mathfrak{B}{(\Sigma,h^\prime)}$, and the cobordism is diffeomorphic to
$[-1,0] \times B{(\Sigma,h)}$ with a $4$-dimensional $2$-handle
attached along $C \subset \Sigma \times \{t\} \subset \{0\} \times
B{(\Sigma,h)}$ with framing $\pf(C)-1$.
\end{propn}
\begin{rmk} \label{R:Lickorish}
This proposition without the contact structures is simply the familiar
observation due to Lickorish~\cite{Lickorish} that surgery on a knot
$C \subset M$ with framing $-1$ relative to a surface $F \supset C$ is
equivalent to splitting $M$ along $F$ and reglueing with a right twist
along $C$. The requirement that $C$ be homologically nontrivial does
not arise if we ignore the contact structures.
\end{rmk}

One should think of the above two propositions as giving two allowable
moves on the set of pairs $(\Sigma,h)$. We will refer to the first
move as ``attaching a $1$-handle'' and to the second move as
``appending a homologically nontrivial right twist''.
\begin{rmk} \label{R:Moves}
A standard fact about Dehn twists is that, if $g\colon\thinspace
\Sigma \to \Sigma$ is any orientation-preserving homeomorphism and $C
\subset \Sigma$ is a simple closed curve, then $g \circ \tau_C \circ
g^{-1} = \tau_{g(C)}$. Now suppose we are given a pair $(\Sigma,h)$,
where $h = h_1 \circ h_2$ for $h_1,h_2 \in \mathcal{M}(\Sigma)$, and
let $D = h_2^{-1}(C)$. Then the operation of replacing $(\Sigma,h)$
with $(\Sigma,h_1 \circ \tau_C \circ h_2)$ is the same as replacing
$(\Sigma,h)$ with $(\Sigma, h \circ \tau_D)$. Thus we can immediately
generalize the second move to allow ``inserting a homologically
nontrivial right twist''. In other words, if $h^\prime = h_1 \circ
\tau_C \circ h_2$, then $\mathfrak{B}{(\Sigma,h)} \prec
\mathfrak{B}{(\Sigma,h^\prime)}$.
\end{rmk}

The next result follows from earlier work by this author~\cite{Gay} on
attaching symplectic $2$-handles along transverse knots; here we will
attach $2$-handles along the binding of an open book
decomposition. Notice that each component $K$ of the binding is given
a framing as a boundary component of a page; call this the
page-framing for the binding and again abbreviate it $\pf(K)$. By a
right twist about a component $C$ of $\partial \Sigma$, we really mean
a right twist about a curve parallel to $C$.
\begin{propn} \label{P:PosTwoHandle}
Let $h^\prime = \delta \circ h^{-1}$, where $\delta$ is the product of
one right twist around each component of $\partial \Sigma$. Then
$(\mathfrak{B}{(\Sigma,h)} \amalg \mathfrak{B}{(\Sigma,h^\prime)}) \prec
\emptyset$. The cobordism is diffeomorphic to $[-1,0] \times
B{(\Sigma,h)}$ with a $4$-dimensional $2$-handle attached to $\{0\}
\times B{(\Sigma,h)}$ along each component $K$ of the binding with
framing $\pf(K)+1$.
\end{propn}
In section~\ref{S:ConstructBlocks}, we will show how this follows from
the results in~\cite{Gay}.

A key point in our construction will depend on a relation among right
twists called the ``chain relation'' by Wajnryb~\cite{Wajnryb}. A
sequence of simple closed curves $\{C_1, \ldots, C_n\}$ in a surface
$\Sigma$ is called a {\em chain} if each $C_i$ intersects each
$C_{i+1}$ transversely in exactly one point.
\begin{propn} [Lemma 21 in~\cite{Wajnryb}] \label{P:Chains}
Let $\{ C_1, \ldots, C_n \}$ be a chain in a surface
$\Sigma$ and let $N$ be a regular neighborhood of $C_1 \cup \ldots
\cup C_n$. If $n$ is even then $\partial N$ has one component. Let
$\tau_i$ be a right twist about $C_i$ and let $\delta$ be a right
twist about $\partial N$. Then the following relation holds:
\[ (\tau_1 \circ \tau_2 \circ \cdots \circ \tau_n)^{2n+2} = \delta \]
\end{propn}
\begin{rmk}
There is a similar relation when $n$ is odd and $\partial N$ has two
components, which could be used here as well, but for simplicity we
only work with the case where $n$ is even.
\end{rmk}
\begin{rmk}
Notice that the regular neighborhood $N$ is actually a surface of
genus $n/2$ with one disk removed. Any surface $\Sigma$ of genus $g$
with one boundary component contains a chain $\{C_1, \ldots C_{2g}\}$
such that $\Sigma$ is a regular neighborhood of $C_1 \cup \ldots \cup
C_{2g}$. Given such a $\Sigma$ and a fixed homologically nontrivial
simple closed curve $C \subset \Sigma$, we may choose the chain so
that $C = C_1$. Furthermore all the curves in such a chain are
necessarily homologically nontrivial.
\end{rmk}

\section{How to put the building blocks together}

\begin{lemma} \label{L:RightsToTheBoundary}
Given a pair $(\Sigma,h)$, where $\Sigma$ is a compact surface with
$\partial \Sigma \neq \emptyset$ and $h \in \mathcal{M}(\Sigma)$,
there exists a pair $(\Sigma_1,h_1)$ obtained from $(\Sigma,h)$ by
adding $1$-handles and inserting homologically nontrivial right twists
(see remark~\ref{R:Moves}) such that:
\begin{itemize}
\item $\Sigma_1$ has one boundary component and,
\item $h_1 = \delta \circ R^{-1}$, where $\delta$ is a right twist
about $\partial \Sigma_1$ and $R$ is a composition of homologically
nontrivial right twists.
\end{itemize}
\end{lemma}
\begin{proof}
By attaching $1$-handles, we can get from $(\Sigma,h)$ to
$(\Sigma_0,h)$ where $\Sigma_0 \supset \Sigma$ has one boundary
component.  By theorem~\ref{T:Dehn}, $h$ is equal to a product of
right-handed and left-handed Dehn twists. First note that each
homologically trivial right or left twist is equal to a product of
homologically nontrivial right or left (respectively) twists, due to
the chain relation (proposition~\ref{P:Chains}). To see this, observe
that any homologically trivial curve $C \subset \Sigma_0$ is the
boundary of a compact subsurface of $\Sigma_0$. Thus we may assume
that all the twists are homologically nontrivial. Let $a$ be the genus
of $\Sigma_0$ and let $\delta_0$ be a right twist about $\partial
\Sigma_0$.

Now suppose we express $h$ as $h = A \circ \tau_C \circ B$ where
$\tau_C$ is a right twist about the homologically nontrival curve $C$
and $A$ and $B$ are arbitrary elements of
$\mathcal{M}(\Sigma_0)$. Then we may find a chain $\{C_1, \ldots,
C_{2a}\}$ in $\Sigma_0$, with $C = C_1$, such that:
\[ \delta_0 = (\tau_1 \circ \ldots \circ \tau_{2a})^{4a+2} \]
Thus by inserting homologically nontrivial right twists we can change
$h$ to $A \circ \delta_0 \circ B = \delta_0 \circ A \circ B$ (since
boundary twists commute with interior twists). Repeating this process
for every right twist in $h$, we can change $h$ to $\delta_0^n \circ
L$ where $L$ is a product of homologically nontrivial left twists and
$n$ is some (possibly large) positive integer. Our task is to reduce
$n$ to $1$. Note that, by inserting more right twists if necessary, we
may assume that $n$ is odd. We write $L = R^{-1}$ where $R$ is the
corresponding product of right twists.

Add more $1$-handles to $\Sigma_0$ to get $\Sigma_1 \supset \Sigma_0$,
with $\partial \Sigma_1$ connected and with the genus $b$ of
$\Sigma_1$ equal to $a n + (n-1)/2$, so that $(4 a + 2)n = 4b +
2$. Fix a particular chain $\{C_1, \ldots, C_{2b}\}$ in $\Sigma_1$
such that $\Sigma_0$ is a regular neighborhood of $\{C_1, \ldots,
C_{2a}\}$ and $\Sigma_1$ is a regular neighborhood of $\{C_1,
\ldots, C_{2b}\}$. Then we have:
\[ \delta_0^n = (\tau_1 \circ \ldots \circ \tau_{2a})^{(4a+2)n} =
(\tau_1 \circ \ldots \circ \tau_{2a})^{4b+2} \]
Thus, by inserting homologically nontrivial right twists we can change
this expression to:
\[ \delta = (\tau_1 \circ \ldots \circ \tau_{2b})^{4b+2} \]
Finally this shows that $\delta_0^n \circ L = \delta_0^n \circ R^{-1}$
may be changed to $\delta \circ R^{-1}$ by inserting right twists.
\end{proof}

\begin{proof}[Proof of theorem~\ref{T:ConcaveFillings}]
Given a contact $3$-manifold $(M,\xi)$, theorem~\ref{T:CSgivesOBD}
tells us that $(M,\xi)$ is supported by an open book decomposition,
with page $\Sigma$ and monodromy $h$. Apply
lemma~\ref{L:RightsToTheBoundary} to get a new pair $(\Sigma_1,h_1)$
with $h_1 = \delta \circ R^{-1}$ as in the lemma. Since we only added
$1$-handles and inserted right twists, propositions~\ref{P:OneHandle}
and~\ref{P:NegTwoHandle} tell us that $(M,\xi) = \mathfrak{B}{(\Sigma,h)}
\prec \mathfrak{B}{(\Sigma_1,h_1)}$.  Consider the pair $(\Sigma_1,\delta
\circ h_1^{-1})$ as in proposition~\ref{P:PosTwoHandle} and notice
that $\delta \circ h_1^{-1} = R$, a composition of homologically
nontrivial right twists. Let $D$ be a disk and $1$ the identity in
$\mathcal{M}(D)$. Then we can get from the pair $(D,1)$ to
$(\Sigma_1,R)$ by attaching $1$-handles and inserting homologically
nontrivial right twists, so that $\emptyset \prec (S^3,\xi_0) \prec
\mathfrak{B}{(\Sigma_1,R)}$ (see remark~\ref{R:SThree}). Thus we have the
following cobordism:
\[ (M,\xi) = (M,\xi) \amalg \emptyset \prec \mathfrak{B}{(\Sigma_1,h_1)}
\amalg \mathfrak{B}{(\Sigma_1,\delta \circ h_1^{-1})} \]
Finally, apply proposition~\ref{P:PosTwoHandle} to conclude that
\[ \mathfrak{B}{(\Sigma_1,h_1)} \amalg \mathfrak{B}{(\Sigma_1,\delta \circ 
h_1^{-1})} \prec \emptyset \]
and put the two cobordisms together to get a concave filling of
$(M,\xi)$. 

\end{proof}

\begin{rmk}
The reader familiar with Lefschetz pencils (see~\cite{GompfStipsicz})
may notice that these building blocks can also be put together to give
a symplectic structure on any topological Lefschetz pencil with
homologically nontrivial vanishing cycles. A Lefschetz pencil is the
result of blowing down a Lefschetz fibration over $S^2$ along $n > 0$
disjoint sections, which must therefore each have self-intersection
$-1$. A Lefschetz fibration over $S^2$ is completely determined by $m$
vanishing cycles $C_1, \ldots, C_m$ in the (closed) fiber surface $F$,
with the property that $\tau_1 \circ \ldots \circ \tau_m = 1 \in
\mathcal{M}(F)$, where $\tau_i$ is a right twist about $C_i$. The
sections correspond to $n$ points $p_1, \ldots, p_n \in F$, disjoint
from the vanishing cycles, such that in fact $\tau_1 \circ \ldots
\circ \tau_m$ is isotopic to the identity via an isotopy fixing
$\{p_1, \ldots, p_n\}$ pointwise. The fact that the sections have
self-intersection $-1$ means that, if we require the isotopy to fix a
disk $D_i$ around each point $p_i$, then $\tau_1 \circ \ldots \circ
\tau_n$ is isotopic to the product of one right twist around each
$\partial D_i$.

Thus the Lefschetz pencil may be built as follows. Remove the interior
of each $D_i$ to get a compact surface $\Sigma$ and let $h = \tau_1
\circ \ldots \circ \tau_n = \delta \in \mathcal{M}(\Sigma)$ (where
$\delta$ is the product of one right twist about each component of
$\partial \Sigma$). First build two copies of a symplectic cobordism
from $\emptyset$ to $\mathfrak{B}{(\Sigma,1)}$ using
remark~\ref{R:SurfaceCrossDisk}. Then attach to one of them a cobordism
from $\mathfrak{B}{(\Sigma,1)}$ to $\mathfrak{B}{(\Sigma,h)}$ built using
proposition~\ref{P:NegTwoHandle} to get a cobordism from $\emptyset$ to
$\mathfrak{B}{(\Sigma,h)}$. Since $h = \delta$,
proposition~\ref{P:PosTwoHandle} gives us a cobordism from
$\mathfrak{B}{(\Sigma,h)} \amalg \mathfrak{B}{(\Sigma,1)}$ to $\emptyset$. Let
$H$ be the union of the $2$-handles used in this last cobordism. These
three cobordisms piece together to make a closed symplectic $4$-manifold
$(X,\omega)$. Since we attached one $2$-handle for each vanishing
cycle to $\Sigma \times D^2$, with the appropriate framings, $X
\setminus H$ is diffeomorphic to the Lefschetz pencil built with this
data with a neighborhood of each section removed. Analysis of the
framings involved in $H$ shows that removing a neighborhood of each
section and replacing them with $H$ is equivalent to blowing down
these sections.

\end{rmk}

\section{Construction of the building blocks} \label{S:ConstructBlocks}

For this section it will be useful to have an explicit model of the
unique contact structure $\xi{(\Sigma,h)}$ supported by an open book
decomposition $(\Sigma,h)$. The following construction is essentially
the construction of Thurston and Winkelnkemper
in~\cite{ThurstonWinkel}, done with a little more care to keep track
of a Reeb vector field. 

Recall that we chose to measure the monodromy of an open book
decomposition using a flow transverse to the pages which had closed
meridinal orbits near the binding. If instead we required closed
longitudinal orbits, realizing the framing $+1$ relative to the page
framing, we would change the monodromy by a single left-handed Dehn
twist along each component of $\partial \Sigma$. Given $(\Sigma,h)$,
let $h^\prime = \delta^{-1} \circ h$, where $\delta$ is a right twist
around each component of $\partial \Sigma$. We will construct
$\mathfrak{B}{(\Sigma,h)}$ using $h^\prime$ as our return map, but
arrange that the return flow orbits are $+1$ longitudes near the
binding.

Given a contact form $\alpha$, let $R(\alpha)$ denote the Reeb vector
field for $\alpha$.

Let $(x,y)$ be coordinates on (each component of) a collar
neighborhood $\nu$ of $\partial \Sigma$, with $x \in (a,b]$, $a>0$, $y
\in S^1$ and $\partial \Sigma = \{x=b\}$. Assume that $h^\prime$ is
the identity on $\nu$. Now choose a $1$-form $\beta$ on $\Sigma$ such
that $d\beta > 0$ and $\beta = x\,dy$ on $\nu$ (this is where we need
$a>0$) and choose a smooth nonincreasing function $f\colon\thinspace
[0,2\pi] \to [0,1]$ which equals $1$ on $[0,2\pi - \epsilon]$ and $0$
on $[2\pi - \epsilon/2,2\pi]$. Let $\alpha = K dt + f(t) \beta +
(1-f(t)) (h^\prime)^* \beta$ on $\Sigma \times [0,2\pi]$, where $t$ is
the $[0,2\pi]$ coordinate and $K$ is some positive constant. Notice
that $\alpha$ descends to a smooth $1$-form on $N = \Sigma \times
[0,2\pi]/\!\sim$, where $(p,2\pi) \sim (h^\prime(p),0)$. For $K$ large
enough, $\alpha \wedge d\alpha > 0$ and $R(\alpha)$ is transverse to
the fibers of the fibration $t\colon\thinspace N \to S^1$. Also notice
that, after choosing $K$, we may enlarge $\Sigma$ to arrange that $x
\in (a,K]$ with $\partial \Sigma = \{x=K\}$ and still have $\beta =
x\,dy$ on $\nu = \{a<x\leq K\}$. Near $\partial N$, $R(\alpha) = (1/K)
\partial_t$.

Now we construct $M$ by Dehn fillings on $N$. For each component of
$\partial N = \partial \Sigma \times S^1$, attach a solid torus $D^2
\times S^1$ to $N \setminus \partial N$ with a map
$\phi\colon\thinspace (D^2 \setminus \{0\}) \times S^1 \to N \setminus
\partial N$ defined as follows: Use polar coordinates
$(r,\mu,\lambda)$ on $D^2 \times S^1$ and coordinates $(x,y,t)$ near
$\partial N$; let $x \circ \phi = (K - r^2)$, $y \circ \phi = -\mu +
\lambda$ and $t \circ \phi = \mu$. This is exactly the right filling
so that the $t$ circles become $+1$ longitudes, so that we are
producing the correct monodromy. Then $\phi^* \alpha = K d\lambda +
r^2 (d\mu - d\lambda)$, which extends across $\{r=0\}$, and $R(\phi^*
\alpha) = (1/K)(\partial_\mu + \partial_\lambda)$. This is tangent to
the binding and satisfies the orientation condition in
definition~\ref{D:Supports}.

Since any two contact structures supported by the same open book
decomposition are isotopic (theorem~\ref{T:OBDGivesCS}), we may always
assume that our contact structures are of the form described above. We
will call this model of the contact structure supported by an open
book decomposition ``the standard model''.

\begin{rmk} \label{R:Legendrian}
A {\em Legendrian} knot is a knot $K$ in a contact $3$-manifold
$(M,\xi)$ which is everywhere tangent to $\xi$. A Legendrian knot $K$
has a canonical framing $\tb(K)$, the ``Thurston-Bennequin'' framing,
given by any vector field in $\xi$ transverse to $K$. If we are given
a homologically nontrivial curve $C \subset \Sigma$, the standard
model may be refined to arrange that $\beta = x\,dy$ in a neighborhood
of $C$, where $x \in (-\epsilon,\epsilon)$, $y \in S^1$ and $C =
\{x=0\}$. Then $C \subset \Sigma \times \{t\} \subset M$ is Legendrian
for $0 < t < 2\pi - \epsilon$, and furthermore $\tb(C) = \pf(C)$.
\end{rmk}

The construction of the $4$-dimensional symplectic handles of
proposition~\ref{P:OneHandle} and proposition~\ref{P:NegTwoHandle} is due to
Weinstein~\cite{Weinstein} and Eliashberg~\cite{EliashStein}; our task
is to show that the contact surgeries associated with these handles
behave well with respect to supporting open book decompositions. This
can be shown by explicit calculations using Weinstein's description of
the handles, but here we give less computational proofs.

\begin{proof}[Proof of proposition~\ref{P:OneHandle}]
Here we essentially ignore Weinstein's construction of a symplectic
$1$-handle and instead build the $1$-handle from scratch.  Let
$(X,\omega)$ be the standard convex filling of
$\mathfrak{B}{(\Sigma,1)}$ and let $(X^\prime,\omega^\prime)$ be the
standard convex filling of $\mathfrak{B}{(\Sigma^\prime,1)}$, as in
remark~\ref{R:SurfaceCrossDisk}. We can clearly construct $(X,\omega)$
and $(X^\prime, \omega^\prime)$ so that $(X,\omega) \subset
(X^\prime,\omega^\prime)$ and so that $H = (X,\omega) \setminus
(X^\prime,\omega^\prime)$ is a $4$-dimensional symplectic
$1$-handle. From the standard model it is clear that there exists a
contactomorphism from a neighborhood $\nu$ of the binding $L \subset
\mathfrak{B}{(\Sigma,h)}$ to a neighborhood $\nu^\prime$ of the
binding $L^\prime \subset \mathfrak{B}{(\Sigma,1)}$. Furthermore we
can arrange that the handle $H$ is attached inside $\nu^\prime$. Thus
$H$ can just as well be attached to a cobordism from
$\mathfrak{B}{(\Sigma,h)}$ to $\mathfrak{B}{(\Sigma,h)}$ to get a
cobordism from $\mathfrak{B}{(\Sigma,h)}$ to
$\mathfrak{B}{(\Sigma^\prime,h^\prime)}$.
\end{proof}

To prove proposition~\ref{P:NegTwoHandle}, we will use Torisu's
characterization in~\cite{Torisu} of the unique contact structure
supported by an open book decomposition.  Suppose that an open book
decomposition of a $3$-manifold $M$ has binding $L$ and fibration
$p\colon\thinspace M \setminus L \to S^1$. Notice that the sets $M_+ =
p^{-1}[0,\pi] \cup L$ and $M_- = p^{-1}[\pi,2\pi] \cup L$ give a
Heegard splitting of $M$. A surface $F \subset M$ is said to be {\em
convex} if there exists a contact vector field $V$ for $\xi$ which is
transverse to $F$. The set $D \subset F$ where $V$ is tangent to $\xi$
is called the {\em dividing set} of $F$. (See~\cite{GirouxConvex} for
background on convexity.) We have the following result, given a
Heegard splitting $(M_+,M_-)$ of $M$ coming from an open book
decomposition as above:
\begin{thm} [Torisu~\cite{Torisu}]
There exists a unique (up to isotopy) contact structure $\xi$ on $M$
such that $\xi$ is tight on both $M_+$ and $M_-$ and such that $F =
\partial M_+ = \partial M_-$ is convex with dividing set equal to the
binding $L$.
\end{thm}
Using the standard model of $(M,\xi) = \mathfrak{B}{(\Sigma,h)}$, one
can show that the splitting surface $F = (\Sigma \times \{0\}) \cup
(\Sigma \times \{\pi\})$ is convex with dividing set $L$. To see that
$\xi$ is tight on both $M_+$ and $M_-$, we note that $(M_+,\xi)$ and
$(M_-,\xi)$ can both be contactomorphically embedded into
$\mathfrak{B}{(\Sigma,1)}$, again using the standard model.
$\mathfrak{B}{(\Sigma,1)}$ is tight because of the existence of the
standard convex filling of $\mathfrak{B}{(\Sigma,1)}$ (see
remark~\ref{R:SurfaceCrossDisk}), and Eliashberg~\cite{EliashFilling}
shows that any contact $3$-manifold with a convex filling is
tight. Putting together the uniqueness in Torisu's theorem with the
uniqueness in theorem~\ref{T:OBDGivesCS}, we get the following
corollary, which in particular shows that Torisu's theorem is actually
equivalent to theorem~\ref{T:OBDGivesCS}.
\begin{cor} \label{C:TorisuSupport}
Given an open book decomposition on $M$, let $M_+$, $M_-$, $F$ and $L$
be as in the preceding paragraphs.  A contact structure $\xi$ on $M$
is supported by this open book decomposition if and only if $\xi$ is
tight on both $M_+$ and $M_-$ and $F$ is convex with dividing set $L$.
\end{cor}

\begin{proof}[Proof of proposition~\ref{P:NegTwoHandle}]

Let $(M,\xi) = \mathfrak{B}{(\Sigma,h)}$. Use the standard model as
described in remark~\ref{R:Legendrian}, so that $C \subset \Sigma
\times \{\pi/2\} \subset (M,\xi)$ is Legendrian, with $\pf(C) =
\tb(C)$. Weinstein~\cite{Weinstein} shows that we can attach a
symplectic handle along an arbitrarily small neighborhood of $C$ with
framing $\tb(C) - 1 = \pf(C) - 1$, so that we get a cobordism from
$(M,\xi)$ to a new contact $3$-manifold $(M^\prime,\xi^\prime)$. As
mentioned earlier, we know that $M^\prime = B{(\Sigma,h^\prime)}$, and
we need to show that $\xi^\prime$ is supported by this open book
decomposition. Split $(M,\xi)$ into the sets $M_+$ and $M_-$ as above;
since the surgery that produces $(M^\prime,\xi^\prime)$ from $(M,\xi)$
is localized near $C \subset M_+$, we get a corresponding Heegard
splitting of $M^\prime$ into two handlebodies $M_+^\prime$ and
$M_-^\prime$, with $(M_-^\prime,\xi^\prime) = (M_-,\xi)$. Since $\xi$
and $\xi^\prime$ also agree on a neighborhood of $F = \partial M_-$,
we still have that the splitting surface $F = F^\prime = \partial
M_+^\prime = \partial M_-^\prime$ is convex with dividing set equal to
the binding $L^\prime$. Thus, by corollary~\ref{C:TorisuSupport}, we
need only show that $(M_+^\prime,\xi^\prime)$ is tight to complete the
proof. Let $\phi$ be the contactomorphic embedding of $(M_+,\xi)$ into
$\mathfrak{B}{(\Sigma,1)}$. We can also attach a symplectic handle
along $\phi(C)$ to the standard convex filling of
$\mathfrak{B}{(\Sigma,1)}$ to get a convex filling of a new contact
$3$-manifold $(M^{\prime\prime},\xi^{\prime\prime})$, which is
therefore tight. Since the contact surgery along $\phi(C)$ must be the
same as the surgery along $C$, we see that $(M_+^\prime,\xi^\prime)$
embeds contactomorphically in
$(M^{\prime\prime},\xi^{\prime\prime})$. Therefore
$(M_+^\prime,\xi^\prime)$ is tight.
\end{proof}

Before proving proposition~\ref{P:PosTwoHandle}, we provide a summary
of the relevant definitions and results from~\cite{Gay}. In that paper
we made a definition very similar to definition~\ref{D:Supports},
except that we worked with structures more general than open book
decompositions. Suppose $(M,\xi)$ is a contact $3$-manifold, $L
\subset M$ is a link and $p\colon\thinspace M \setminus L \to S^1$ is
a fibration.
\begin{defn} \label{D:FakeSupport}
The pair $(L,p)$ is a nicely fibered link supporting $\xi$ if there
exists a Reeb field $V$ for $\xi$ with $dp(V) >0$ and polar
coordinates $(r,\mu,\lambda)$ near each component of $L$ such that the
following conditions are satisfied on each coordinate neighborhood of $L$:
\begin{itemize}
\item $\xi = \ker( f(r) d\lambda + g(r) d\mu)$ for some functions $f(r)$
and $g(r)$.
\item $V = A \partial_\mu + B \partial_\lambda$ and $p = C \mu
+ D \lambda$, where $A,B,C,D$ are constant near each component of $L$
and $B$ and $C$ are positive.
\end{itemize}
\end{defn}
Note that this definition implies that $V$ is tangent to $L$. We may
think of $(L,p)$ as a ``fake open book decomposition'', since the
boundary of a ``page'' may multiply cover the ``binding''. There
appears to be much more stringent control of $\xi$ near $L$ than in
definition~\ref{D:Supports} and the orientation condition (that $B$
and $C$ are positive) looks different. However the following lemma can
be proved by direct computation using the standard model.
\begin{lemma}
If the fibers of $p$ meet $L$ as longitudes so that we have an honest
open book decomposition of $M$, then the condition that $\xi$ is
supported by the open book decomposition is equivalent, after an
isotopy, to the condition that $(L,p)$ is a nicely fibered link
supporting $\xi$.
\end{lemma}

Notice that we can compare framings of $L$ to the fibration $p$. In
terms of polar coordinates $(r,\mu,\lambda)$ near $L$, such a
fibration $p$ determines a ``slope'' $d\mu/d\lambda = s_p =
-D/C \in \mathbb{R} \cup \{\infty\}$, while a framing $f$
determines a family of parallel longitudes with ``slope''
$d\mu/d\lambda = s_f \in \mathbb{Z}$.
\begin{defn}
We say that the framing $f$ is positive with respect to $p$ if, on
each component of $L$, $s_f > s_p$.
\end{defn}
\begin{rmk}
When $p$ gives an honest open book decomposition, $f$ is positive with
respect to $p$ if, on each component $K$ of $L$, $f = \pf(K) + k$ for
some positive integer $k$.
\end{rmk}

Following is the main result we need from~\cite{Gay}. This is
essentially theorem~1.2 and addendum~5.1 in~\cite{Gay}, but we also
include some points that are made in their proofs.
\begin{thm} \label{T:Thesis}
Suppose that $(L,p)$ is a nicely fibered link in $M$ supporting a
contact structure $\xi$ and that $f$ is a framing of $L$ which is
positive with respect to $p$. Let $X$ be $[-1,0] \times M$ with a
$2$-handle attached along each component of $L$ in $\{0\} \times M$
with framing $f$. Let $M = \{-1\} \times M$ and let $M_1$ be the other
component of $\partial X$, both oriented opposite to the boundary
orientation. Let $L_1 \subset M_1$ be the union of the ascending
spheres of the $2$-handles. Then $X$ supports a symplectic structure
$\omega$ and a Liouville vector field $V$ defined near $\partial X$
pointing in along both $M$ and $M_1$ and inducing the given contact
structure $\xi$ on $M$. Let $\alpha$ and $\alpha_1$ be the
restrictions of $\imath_V \omega$ to $M$ and $M_1$, respectively. Then
we have the following relation between $\alpha$ and $\alpha_1$:
\begin{itemize}
\item There exist polar coordinates $(r,\mu,\lambda)$ on a
neighborhood $\nu$ of $L$ as in definition~\ref{D:FakeSupport}, a
closed tubular neighborhood $\tau = \{r \leq \epsilon\} \subset \nu$
of $L$, a positive constant $k$, a function $h\colon\thinspace M
\setminus \tau \to [1,\infty)$ which is $1$ outside $\nu$ and a
function of $r$ inside $\nu$, and an orientation-reversing
diffeomorphism $\phi\colon\thinspace M \setminus \tau \to M_1
\setminus L_1$ such that: $\phi^*(\alpha_1) = k\,dp - h \alpha$
\item The pair $(L_1, p_1 = (\phi^{-1})^* p)$ is a nicely fibered link
supporting $\xi_1 = \ker \alpha_1$.
\end{itemize}
\end{thm}
Here is the idea of the relationship between $(M,\xi)$ and
$(M_1,\xi_1)$: First we perform the topological surgery by removing
the interior of $\tau$ from $M$ and then collapsing each component of
$\partial \tau$ to a circle, to get the new link $L_1 \subset M_1$. The
theorem makes a judicious choice of $k$ and $h$, so as to arrange that
$k\,dp - e^h \alpha$ is a negative contact form on $M \setminus \nu$
and furthermore extends across $L_1$ to a contact form $\alpha_1$ on
$M_1$. Then we reverse the orientation to treat $\xi_1 = \ker
\alpha_1$ as a positive contact structure on $M_1$.

\begin{proof}[Proof of proposition~\ref{P:PosTwoHandle}]
Let $(M,\xi) = \mathfrak{B}{(\Sigma,h)}$, with binding $L$ and
fibration \linebreak $p\colon\thinspace M \setminus L \to
S^1$. Proposition~\ref{P:PosTwoHandle} asks us to attach a handle
along each component $K$ of $L$ with framing $\pf(K)+1$.
Theorem~\ref{T:Thesis} tells us that we get a symplectic cobordism
from $(M,\xi) \amalg (M_1,\xi_1)$ to $\emptyset$. We must show that
$(M_1,\xi_1)$ is supported by the open book decomposition
$(\Sigma,\delta \circ h^{-1})$. Because the framing of the surgery is
$\pf(K)+1$ for each component $K$, the fibers of $p_1$ still meet
$L_1$ as longitudes, so that we do have an honest open book
decomposition of $M_1$ which supports $\xi_1$ and with pages
diffeomorphic to $\Sigma$. To compute the monodromy, note that a
meridian before the surgery becomes a longitude with framing
$\pf(L_1)-1$, so that our original monodromy $h$ is now measuring the
monodromy via a flow which is longitudinal near $L_1$. Correcting this
flow to be meridinal changes the monodromy by a left twist on each
boundary component, so that, properly measured, the monodromy for
$-M_1$ is $\delta^{-1} \circ h$. We should reverse the flow to get the
monodromy for $M_1$, which is thus $\delta \circ h^{-1}$.
\end{proof}

\bibliography{CPSconcave}

\bibliographystyle{abbrv}

%\end{spacing}{2}

\end{document}